# A new concentration result for regularized risk minimizers

**Ingo Steinwart,**[1] **Don Hush**[1] **and Clint Scovel**[1]

*Los Alamos National Laboratory*

**Abstract:** We establish a new concentration result for regularized risk minimizers which is similar to an oracle inequality. Applying this inequality to regularized least squares minimizers like least squares support vector machines, we show that these algorithms learn with (almost) the optimal rate in some specific situations. In addition, for regression our results suggest that using the loss function $L_\alpha(y,t) = |y-t|^\alpha$ with $\alpha$ near 1 may often be preferable to the usual choice of $\alpha = 2$.

## 1. Introduction

The theoretical understanding of support vector machines (SVMs) and related kernel-based methods has been substantially improved in recent years. Based on Talagrand's concentration inequality and local Rademacher averages it has recently been shown that SVMs for classification can learn with rates up to $\frac{1}{n}$ under somewhat realistic assumptions on the data-generating distribution (see [12] and the related work [3]). However, the currently available technique, namely the so-called "shrinking technique" in [12], for establishing such rates requires choosing the *entire* regularization sequence *a-priori*. Unfortunately, the optimal regularization sequences usually depend on some features of the data-generating distribution typically unknown in practice, and consequently the results derived by the shrinking technique have some serious drawbacks.

In this work we replace the shrinking technique by a localization argument similar to the localization argument used in conjunction with local Rademacher averages. The key observation for this new localization argument is that regularized risk minimizers control the size of the norm in the regularization term by their (excess) risk in a non-trivial manner (see Lemma 4.1 for details). As a consequence of this observation, we can not only localize with respect to small variances but also with respect to small maximum norms.

Using the above (double) localization we obtain oracle-type inequalities for a large class of regularized risk minimizers including support vector machines, and regularization networks. For the former we can easily reproduce rates established in [12, 13], while for the latter we show some minmax rates in specific situations and provide results indicating that using the loss function $L_\alpha(y,t) = |y-t|^\alpha$ with $\alpha$ near 1 to estimate the regression function may be more robust to both outliers and the choice of regularization parameter than the usual choice $\alpha = 2$.

## 2. An oracle inequality for regularized risk minimizers

Throughout this work we assume that $X$ is compact metric space, $Y \subset [-1,1]$ is compact, $P$ is a Borel probability measure on $X \times Y$, and $H$ is a RKHS of continuous functions over $X$ with closed unit ball $B_H$. It is well-known that $H$ can then be continuously embedded into the space of continuous functions $C(X)$ equipped with the usual maximum-norm $\|.\|_\infty$. In order to avoid constants we always assume that this embedding has norm 1, i.e. $\|.\|_\infty \leq \|.\|_H$.

[1]Modeling, Algorithms and Informatics Group, CCS-3, Los Alamos National Laboratory, Los Alamos, New Mexico 87545, USA, e-mail: ingo@lanl.gov; dhush@lanl.gov; jcs@lanl.gov

Furthermore, $L : Y \times \mathbb{R} \to [0, \infty)$ always denotes a continuous function which is convex in the second variable. In the following we are particularly interested in functions $L$ that satisfy the growth assumptions introduced in [6]:

(1) $\quad \sup_{y \in Y} L(y, t) \leq 1 + |t|^\alpha \qquad \text{and} \qquad \sup_{y \in Y} \big|L_{|Y \times [-t,t]}(y, .)\big|_1 \leq c_L \, t^{\alpha - 1}$

for some constants $\alpha \in [1, 2]$, $c_L > 0$, and all $t \in \mathbb{R}$, where $|h|_1$ denotes the Lipschitz constant of a function $h$. The functions $L$ will serve as loss functions and consequently let us recall the associated $L$-risk

$$\mathcal{R}_{L,P}(f) = \mathbb{E}_{(x,y) \sim P} L(y, f(x)),$$

where $f : X \to \mathbb{R}$ is a measurable function. Note that (1) immediately gives $\mathcal{R}_{L,P}(0) \leq 1$. Furthermore, the minimal $L$-risk is denoted by $\mathcal{R}^*_{L,P}$, i.e.

$$\mathcal{R}^*_{L,P} = \inf\{\mathcal{R}_{L,P}(f) \,|\, f : X \to \mathbb{R} \text{ measurable}\},$$

and a function attaining this infimum is denoted by $f^*_{L,P}$.

The learning schemes we are interested in are based on an optimization problem of the form

$$f_{P,\lambda} := \arg\min_{f \in H} \left(\lambda \|f\|_H^2 + \mathcal{R}_{L,P}(f)\right),$$

where $\lambda > 0$. Note that if we identify a training set $T = ((x_1, y_1), \ldots, (x_n, y_n)) \in (X \times Y)^n$ with its empirical measure, then $f_{T,\lambda}$ denotes the empirical estimators of the above learning scheme. Obviously, support vector machines (see e.g. [5]) and regularization networks (see e.g. [8]) are both learning algorithms which fall into the above category.

One way to describe the approximation error of these learning schemes is the *approximation error function*

$$a(\lambda) := \lambda \|f_{P,\lambda}\|^2 + \mathcal{R}_{L,P}(f_{P,\lambda}) - \mathcal{R}^*_{L,P}, \qquad \lambda > 0,$$

which we discussed in some detail in [13]. Furthermore in order to deal with the complexity of the used RKHSs let us recall that for a subset $A \subset E$ of a Banach space $E$ the *covering numbers* are defined by

$$\mathcal{N}(A, \varepsilon, E) := \min\Big\{n \geq 1 : \exists x_1, \ldots, x_n \in E \text{ with } A \subset \bigcup_{i=1}^n (x_i + \varepsilon B_E)\Big\}, \qquad \varepsilon > 0,$$

where $B_E$ denotes the closed unit ball of $E$. Given a finite sequence $T = (z_1, \ldots, z_n) \in Z^n$ we are particularly interested in the Banach space $L_2(T)$ which consists of all equivalence classes of functions $f : Z \to \mathbb{R}$ and which is equipped with the norm

(2) $\quad \|f\|_{L_2(T)} := \Big(\frac{1}{n} \sum_{i=1}^n |f(z_i)|^2\Big)^{\frac{1}{2}}.$

In other words, $L_2(T)$ is a $L_2$-space with respect to the empirical measure of $(z_1, \ldots, z_n)$. Furthermore, if $T$ is of the form $T = ((x_1, y_1), \ldots, (x_n, y_n))$, and $T_X := (x_1, \ldots, x_n)$, then the space $L_2(T_X)$ has the obvious meaning. In addition to the convention $0^0 := 1$ we utilize the following

$$(3) \qquad a^\infty := \begin{cases} 0 & \text{if } 0 \leq a < 1, \\ 1 & \text{if } a = 1, \\ \infty & \text{if } a > 1. \end{cases}$$

Now we can state the main result of this paper:

**Theorem 2.1.** *Let $H$ be a RKHS of a continuous kernel over $X$ with $\|.\|_\infty \leq \|.\|_H$. Assume that there are constants $a \geq 1$ and $0 < p < 2$ such that for all $\delta > 0$ we have*

$$(4) \qquad \sup_{T \in Z^n} \log \mathcal{N}(B_H, \delta, L_2(T_X)) \leq a\delta^{-p}.$$

*Let $L : Y \times \mathbb{R} \to [0, \infty)$ be a continuous function which is convex in its second variable and satisfies (1). Furthermore, let $P$ be a distribution on $X \times Y$ such that $f^*_{L,P}$ exists. Moreover, suppose that for all $0 < \lambda \leq 1$ and all $f \in \lambda^{-\frac{1}{2}} B_H$ we have*

$$(5) \qquad \mathbb{E}_P \big(L \circ f - L \circ f^*_{L,P}\big)^2 \leq c \big(\|f\|_\infty + 1\big)^v \big(\mathbb{E}_P L \circ f - L \circ f^*_{L,P}\big)^\vartheta$$

*for some constants $c \geq 1$, $\vartheta \in (0, 1]$, and $v \in [0, 2]$. Then there exists a constant $K \geq 1$ such that for all $0 < \lambda \leq 1$, $\varepsilon > 0$, $x \geq 1$ satisfying*

$$\varepsilon \geq \max\bigg\{ a(\lambda) + \lambda, \bigg(\frac{Ka}{\lambda^{\frac{2\alpha p + v(2-p)}{4}} n}\bigg)^{\frac{4}{8-2\alpha p - (v+2\vartheta)(2-p)}}, \bigg(\frac{Ka}{\lambda^{\frac{\alpha(2+p)}{4}} n}\bigg)^{\frac{4}{(2+p)(2-\alpha)}},$$

$$\bigg(\frac{Kx}{\lambda^{\frac{v}{2}} n}\bigg)^{\frac{2}{4-v-2\vartheta}}, \bigg(\frac{Kx}{\lambda^{\frac{\alpha}{2}} n}\bigg)^{\frac{2}{2-\alpha}}\bigg\},$$

*we have*

$$\Pr^*\Big(T \in Z^n : \mathcal{R}_{L,P}(f_{T,\lambda}) - \mathcal{R}^*_{L,P} < a(\lambda) + \varepsilon\Big) \geq 1 - e^{-x}$$

*where $\Pr^*$ denotes the outer probability.*

Theorem 2.1 is proved in Section 4. Now we proceed to illustrate its utility with some applications.

**Example 2.2** (Least square regression with Sobolev spaces)**.** Let us consider the least squares loss function which is defined by $L(y, t) = (y - t)^2$. Furthermore, let us assume that $H$ contains the regression function $x \mapsto \mathbb{E}(y|x)$ and satisfies the complexity exponent condition (4). In addition let $(\lambda_n)$ be a strictly positive null-sequence with $\lambda_n^{1+p/2} n \to \infty$. Then in Section 5 we show that our learning rate is of the form $\lambda_n$. In particular, if $H$ is a Sobolev space of order $m$ on some suitable $X \subset \mathbb{R}^d$, $m > d/2$, then we have $p = d/m$, and consequently, for $\lambda_n := n^{-\frac{2m}{2m+d}} \log n$ our rate becomes $n^{-\frac{2m}{2m+d}} \log n$. This equals the optimal rate $n^{-\frac{2m}{2m+d}}$ up to a logarithmic factor (see e.g. [7] and the references therein).

**Example 2.3** (Comparison of different loss functions used for regression)**.** Consider again regression with the squared loss function $L_2(y, t) = (y - t)^2$ defining performance but use the loss function $L_\alpha(y, t) = |y - t|^\alpha$ with $1 \leq \alpha \leq 2$ to determine the estimate $f_{T,\lambda}$. Suppose that $H$ contains the regression function $x \mapsto \mathbb{E}(y|x)$, and satisfies the complexity exponent condition (4). In Section 5 we begin by using the oracle inequality of Theorem 2.1 to bound the excess $L_\alpha$-risk $\mathcal{R}_{L_\alpha,P}(f_{T,\lambda}) - \mathcal{R}^*_{L_\alpha,P}$. When $\alpha = 2$ we produce the results of Example 2.2. When $1 < \alpha < 2$ we set $\lambda = n^{-\kappa}$ with $\kappa > 0$ and observe that when $\kappa \leq \frac{2}{2+p}$ we obtain the rate $n^{-\kappa}$ independently of the value of $\alpha$ and when $\kappa > \frac{2}{2+p}$ we obtain the rate $n^{-\frac{2}{2+p}+(\kappa-\frac{2}{2+p})\frac{\alpha}{2-\alpha}}$. We conclude that the $\kappa$-optimal learning rate for the $L_\alpha$ risk is $n^{-\frac{2}{2+p}}$ and is achieved when $\kappa = \frac{2}{2+p}$. Now suppose that the conditional distributions $P(y|x)$ are symmetric. These results are then combined with a calibration inequality

$$\mathcal{R}_{L_2,P}(f_{T,\lambda}) - \mathcal{R}^*_{L_2,P} \leq \Psi(\mathcal{R}_{L_\alpha,P}(f_{T,\lambda}) - \mathcal{R}^*_{L_\alpha,P})$$

derived from [11] to obtain bounds on $\mathcal{R}_{L_2,P}(f_{T,\lambda}) - \mathcal{R}^*_{L_2,P}$ in terms of $1 < \alpha < 2$. We observe that when $\kappa \leq \frac{2}{2+p}$ we obtain the rate $n^{-\kappa}$ independently of the value of $\alpha$ and when $\kappa > \frac{2}{2+p}$ we obtain the rate $n^{-\frac{2}{2+p}+(\kappa-\frac{2}{2+p})\frac{2}{2-\alpha}}$. We conclude that the $\kappa$-optimal learning rate for the $L_2$ risk also is $n^{-\frac{2}{2+p}}$ and is achieved when $\kappa = \frac{2}{2+p}$. It is important to observe that the rate for fixed $\kappa$ gets worse as $\alpha$ increases towards 2 and in particular that we have no rates when $2 - (\kappa - \frac{2}{2+p})(2+p) \leq \alpha \leq 2$. When $\alpha = 1$ [11, Example 3.25] shows how, even though the loss function is not strictly convex, we can obtain a calibration inequality in terms of assumptions concerning the concentration about the mean. Consequently with extra assumptions regarding concentration about the mean we can apply these methods, but do not carry out such calculations here since they they are out of the scope of this paper. Moreover, since $\alpha = 1$ is considered more robust to outliers than $\alpha = 2$, these results suggest that setting $\alpha$ near 1 has some substantial advantages to the usual choice $\alpha = 2$. However, to make such a claim more precise will require considering whether and in which sense the assumptions of symmetry and boundedness have been violated. Finally, let us now consider when $H$ is a Sobolev space as in Example 2.2. Then it is clear that we obtain the same optimal rates for all values of $1 < \alpha \leq 2$, although for $\alpha$ near 1 we should concern ourselves with the arising constants.

**Example 2.4** (Hinge loss classification). Let $Y := \{-1, 1\}$, $L$ be defined by $L(y,t) := \max\{0, 1 - yt\}$, $y \in Y$, $t \in \mathbb{R}$, and $P$ be a distribution with Tsybakov noise exponent $q \in [0, \infty]$ in the sense of [12, 13] (see also [2]). When $q > 0$, it follows from [12, Lemma 6.6] that the assumption (5) is satisfied with $\alpha = 1$, $v = \frac{q+2}{q+1}$, $\vartheta = \frac{q}{q+1}$ and $c = \|(2\eta - 1)^{-1}\|_{q,\infty} + 2$. Moreover it is simple to show the same is true when $q = 0$ but with $c = 5$. Hence the condition on $\varepsilon$ becomes

$$\varepsilon \geq \max\left\{a(\lambda) + \lambda, \frac{K}{\lambda}\left(\frac{a}{n}\right)^{\frac{4(q+1)}{2q+pq+4}}, \frac{K}{\lambda}\left(\frac{a}{n}\right)^{\frac{4}{2+p}}, \frac{K}{\lambda}\left(\frac{x}{n}\right)^{\frac{2(q+1)}{q+2}}, \frac{K}{\lambda}\left(\frac{x}{n}\right)^2\right\}.$$

Some easy estimates then show that this reduces to

$$\varepsilon \geq a(\lambda) + \lambda + Kx^2\lambda^{-1}\left(\frac{a}{n}\right)^{\frac{4(q+1)}{2q+pq+4}},$$

where $K \geq 1$ is a suitable constant and $a$ and $n$ are assumed to satisfy $n \geq a \geq 1$. From this we immediately obtain the rates established in [12, Thm. 2.8] and [13, Thm. 1].

## 3. A concentration result for ERM schemes

The proof of our main result Theorem 2.1 is based on a refinement of standard local Rademacher average techniques. Since this refinement may be of its own interest we separate its presentation from the proof of 2.1.

Let us begin by introducing some notations. To this end let $\mathcal{F}$ be a class of bounded measurable functions from $Z$ to $\mathbb{R}$. In order to avoid measurability considerations we always assume that *$\mathcal{F}$ is separable with respect to $\|.\|_\infty$*. Given a probability measure $P$ on $Z$ we define the modulus of continuity of $\mathcal{F}$ by

$$\omega_{P,n}(\mathcal{F}, \varepsilon) := \mathbb{E}_{T \sim P^n} \left( \sup_{\substack{f \in \mathcal{F}, \\ \mathbb{E}_P f \leq \varepsilon}} |\mathbb{E}_P f - \mathbb{E}_T f| \right),$$

where we emphasize that the supremum is, as a function from $Z^n$ to $\mathbb{R}$, measurable by the separability assumption on $\mathcal{F}$. In addition note that the supremum is taken over all $f \in \mathcal{F}$ with $\mathbb{E}_P f \leq \varepsilon$, whereas usually the supremum is taken over all $f \in \mathcal{F}$ with $\mathbb{E}_P f^2 \leq \varepsilon$.

We also need some notations related to ERM-type algorithms: we call $C : \mathcal{F} \times Z \to [0,\infty)$ a *cost function* if $C \circ f := C(f,.)$ is measurable for all $f \in \mathcal{F}$. Given a probability measure $P$ on $Z$ we denote by $f_{P,\mathcal{F}} \in \mathcal{F}$ a minimizer of

$$f \mapsto \mathcal{R}_{C,P}(f) := \mathbb{E}_{z \sim P} C(f,z).$$

Moreover, if $P$ is an empirical measure with respect to $T \in Z^n$ we write $f_{T,\mathcal{F}}$ and $\mathcal{R}_{C,T}(.)$ as usual. For simplicity, we assume throughout this section that $f_{P,\mathcal{F}}$ and $f_{T,\mathcal{F}}$ do exist. Furthermore, although there may be multiple solutions we use a single symbol for them whenever no confusion regarding the non-uniqueness of this symbol can be expected. An algorithm that produces solutions $f_{T,\mathcal{F}}$ is called an *empirical $C$-risk minimizer*. Moreover, if $\mathcal{F}$ is convex, we say that $C$ is convex if $C(.,z)$ is convex for all $z \in Z$. Finally, $C$ is called *line-continuous* if for all $z \in Z$ and all $f, \hat{f} \in \mathcal{F}$ the function $t \mapsto C(tf + (1-t)\hat{f}, z)$ is continuous on $[0,1]$. If $\mathcal{F}$ is a vector space then every convex $C$ is line-continuous. Now we can formulate the main result of this section:

**Theorem 3.1.** *Let $\mathcal{F}$ be a convex set of bounded measurable functions from $Z$ to $\mathbb{R}$, $C : \mathcal{F} \times Z \to [0,\infty)$ be a convex, line-continuous cost function, and $P$ be a probability measure on $Z$. Assume that*

$$\mathcal{G} := \{C \circ f - C \circ f_{P,\mathcal{F}} \; : \; f \in \mathcal{F}\}$$

*is separable with respect to $\|.\|_\infty$. Furthermore assume that there exist constants $b, B \geq 0$, $\beta \in [0,1]$, and $w, W \geq 0$, $\nu \in [0,2]$, $\vartheta \in [0,2)$, such that*

(6) $$\|g\|_\infty \leq b \left(\mathbb{E}_P g\right)^\beta + B$$

*and*

(7) $$\mathbb{E}_P g^2 \leq \left(b(\mathbb{E}_P g)^\beta + B\right)^\nu \left(w(\mathbb{E}_P g)^\vartheta + W\right)$$

*for all $g \in \mathcal{G}$. Then for $n \geq 1$, $x \geq 1$ and $\varepsilon > 0$ satisfying*

$$\varepsilon \geq 3\omega_{P,n}(\mathcal{G}, \varepsilon) + \sqrt{\frac{2x(b\varepsilon^\beta + B)^\nu(w\varepsilon^\vartheta + W)}{n}} + \frac{2x(b\varepsilon^\beta + B)}{n}$$

*we have*

$$\Pr^* \left( T \in Z^n : \mathcal{R}_{C,P}(f_{T,\mathcal{F}}) < \mathcal{R}_{C,P}(f_{P,\mathcal{F}}) + \varepsilon \right) \geq 1 - e^{-x}.$$

In order to prove Theorem 3.1 let us first recall Talagrand's concentration inequality (see [14]). The following version of this inequality is derived from Bousquet's result in [4] using a little trick presented in [1, Lem. 2.5]:

**Theorem 3.2.** *Let $P$ be a probability measure on $Z$ and $\mathcal{H}$ be a set of bounded measurable functions from $Z$ to $\mathbb{R}$ which is separable with respect to $\|.\|_\infty$ and satisfies $\mathbb{E}_P h = 0$ for all $h \in \mathcal{H}$. Furthermore, let $M > 0$ and $\tau \geq 0$ be constants with $\|h\|_\infty \leq M$ and $\mathbb{E}_P h^2 \leq \tau$ for all $h \in \mathcal{H}$. Then for all $x \geq 1$ and all $n \geq 1$ we have*

$$P^n\left(T \in Z^n : \sup_{h \in \mathcal{H}} \mathbb{E}_T h > 3\mathbb{E}_{T' \sim P^n} \sup_{h \in \mathcal{H}} \mathbb{E}_{T'} h + \sqrt{\frac{2x\tau}{n}} + \frac{Mx}{n}\right) \leq e^{-x}.$$

This concentration inequality is used to prove the following lemma which is a generalized version of Lemma 13 in [2] and Lemma 5.4 in [12]:

**Lemma 3.3.** *Let $P$ be a probability measure on $Z$ and $\mathcal{G}$ be a set of bounded measurable functions from $Z$ to $\mathbb{R}$ which is separable with respect to $\|.\|_\infty$. Let us assume that $\mathcal{G}$ satisfies (6) and (7), and that there is a constant $a \in [0, 1)$ such that for all $T \in Z^n$, $\varepsilon > 0$ for which there is a $g \in \mathcal{G}$ with*

$$\mathbb{E}_T g \leq a\varepsilon \qquad and \qquad \mathbb{E}_P g \geq \varepsilon$$

*there is also an element $g^* \in \mathcal{G}$ with*

$$\mathbb{E}_T g^* \leq a\varepsilon \qquad and \qquad \mathbb{E}_P g^* = \varepsilon.$$

*Then for all $n \geq 1$, $x \geq 1$, and all $\varepsilon > 0$ satisfying*

$$(1-a)\varepsilon \geq 3\omega_{P,n}(\mathcal{G}, \varepsilon) + \sqrt{\frac{2x(b\varepsilon^\beta + B)^\nu(w\varepsilon^\vartheta + W)}{n}} + \frac{2x(b\varepsilon^\beta + B)}{n}$$

*we have*

$$\Pr{}^*\left(T \in Z^n : \text{ for all } g \in \mathcal{G} \text{ with } \mathbb{E}_T g \leq a\varepsilon \text{ we have } \mathbb{E}_P g < \varepsilon\right) \geq 1 - e^{-x}.$$

*Proof.* We define $\mathcal{H} := \{\mathbb{E}_P g - g : g \in \mathcal{G}, \mathbb{E}_P g = \varepsilon\}$. Obviously, for all $h \in \mathcal{H}$ we have $\mathbb{E}_P h = 0$ and

$$\|h\|_\infty \leq 2b\varepsilon^\beta + 2B =: M,$$
$$\mathbb{E}_P h^2 \leq \mathbb{E}_P g^2 \leq (b\varepsilon^\beta + B)^\nu(w\varepsilon^\vartheta + W) =: \tau.$$

Moreover, it is also easy to verify that $\mathcal{H}$ is separable with respect to $\|.\|_\infty$. As in the proof of Lemma 5.4 in [12] our assumption on $\mathcal{G}$ now yields

$$\Pr{}^*\left(T \in Z^n : \exists g \in \mathcal{G} \text{ with } \mathbb{E}_T g \leq a\varepsilon \text{ and } \mathbb{E}_P g \geq \varepsilon\right)$$
$$\leq \Pr{}^*\left(T \in Z^n : \exists g \in \mathcal{G} \text{ with } \mathbb{E}_T g \leq a\varepsilon \text{ and } \mathbb{E}_P g = \varepsilon\right)$$
$$= \Pr{}^*\left(T \in Z^n : \exists g \in \mathcal{G} \text{ with } \mathbb{E}_P g - \mathbb{E}_T g \geq (1-a)\varepsilon \text{ and } \mathbb{E}_P g = \varepsilon\right)$$
$$\leq P^n\left(T \in Z^n : \sup_{\substack{g \in \mathcal{G} \\ \mathbb{E}_P g = \varepsilon}} (\mathbb{E}_P g - \mathbb{E}_T g) \geq (1-a)\varepsilon\right)$$
$$= P^n\left(T \in Z^n : \sup_{h \in \mathcal{H}} \mathbb{E}_T h \geq (1-a)\varepsilon\right).$$

In order to bound the last probability we will apply Theorem 3.2. To this end observe

$$3\mathbb{E}_{T'\sim P^n}\sup_{h\in\mathcal{H}}\mathbb{E}_{T'}h + \sqrt{\frac{2x\tau}{n}} + \frac{Mx}{n} \leq (1-a)\varepsilon\,,$$

and consequently applying Theorem 3.2 yields

$$\Pr{}^*\bigl(T\in Z^n : \exists g\in\mathcal{G} \text{ with } \mathbb{E}_T g \leq a\varepsilon \text{ and } \mathbb{E}_P g \geq \varepsilon\bigr) \leq e^{-x}. \qquad \square$$

With the help of the above lemma we can now prove Theorem 3.1:

*Proof of Theorem 3.1.* For $a := 0$ we will apply Lemma 3.3 to the class $\mathcal{G}$. To this end it obviously suffices to show the richness condition on $\mathcal{G}$ of Lemma 3.3: let $f \in \mathcal{F}$ satisfy

$$\mathbb{E}_T(C\circ f - C\circ f_{P,\mathcal{F}}) \leq 0 \qquad \text{and} \qquad \mathbb{E}_P(C\circ f - C\circ f_{P,\mathcal{F}}) \geq \varepsilon\,.$$

For $t \in [0,1]$ we define $f_t := tf + (1-t)f_{P,\mathcal{F}}$. Since $\mathcal{F}$ is convex we have $f_t \in \mathcal{F}$ for all $t \in [0,1]$. By the line-continuity of $C$ and Lebesgue's theorem we find that the map $h : t \mapsto \mathbb{E}_P(C\circ f_t - C\circ f_{P,\mathcal{F}})$ is continuous for $t \in [0,1]$. Since $h(0) = 0$ and $h(1) \geq \varepsilon$ there is a $t \in (0,1]$ with

$$\mathbb{E}_P(C\circ f_t - C\circ f_{P,\mathcal{F}}) \;=\; h(t) \;=\; \varepsilon$$

by the intermediate value theorem. Moreover, for this $t$ the convexity of $C$ gives

$$\mathbb{E}_T(C\circ f_t - C\circ f_{P,\mathcal{F}}) \;\leq\; \mathbb{E}_T\Bigl(tC\circ f + (1-t)C\circ f_{P,\mathcal{F}} - C\circ f_{P,\mathcal{F}}\Bigr) \;\leq\; 0\,.$$

Now, let $\varepsilon > 0$ satisfy the assumption of the theorem. Then $\varepsilon$ also satisfies the assumptions of Lemma 3.3, and hence we find that with probability at least $1-e^{-x}$ every $f \in \mathcal{F}$ with $\mathbb{E}_T(C\circ f - C\circ f_{P,\mathcal{F}}) \leq 0$ satisfies $\mathbb{E}_P(C\circ f - C\circ f_{P,\mathcal{F}}) < \varepsilon$. Since we always have

$$\mathbb{E}_T\bigl(C\circ f_{T,\mathcal{F}} - C\circ f_{P,\mathcal{F}}\bigr) \;\leq\; 0$$

we obtain the assertion. $\qquad\square$

## 4. Proof of the main result

In order to prove our oracle-type inequality we will apply Theorem 3.1. To this end we define the regularized cost function $C_\lambda$ by

$$C_\lambda(x,y,f) := \lambda\|f\|_H^2 + L(y,f(x))\,, \qquad x\in X,\, y\in Y,\, f\in H,$$

and the induced cost class

$$\mathcal{G}(\lambda) := \bigl\{C_\lambda\circ f - C_\lambda\circ f_{P,\lambda} : f \in \lambda^{-1/2}B_H\bigr\}\,, \qquad \lambda > 0.$$

Obviously, the $C_\lambda$-risk minimizer produces the functions $f_{P,\lambda}$ and $f_{T,\lambda}$. Note that $\mathcal{R}_{L,P}(0) \leq 1$ implies $f_{P,\lambda} \in \lambda^{-1/2}B_H$ for all distributions $P$ on $X\times Y$, and hence the latter in particular holds for the empirical solutions $f_{T,\lambda}$. However, it was already observed in [12] that, depending on the approximation error function, sharper bounds for $\|f_{T,\lambda}\|$ are possible with high probability. In order to establish such sharper bounds we employed a "shrinking technique" in [12] which is rather complicated. The key idea of this paper is to replace the shrinking technique by a localization argument based on (6). Consequently, let us first show that regularized risk minimizers always satisfy the supremum bound (6):

**Lemma 4.1.** *Let $0 < \lambda \leq 1$, and suppose that $g \in \mathcal{G}(\lambda)$. Then for any $f \in \lambda^{-\frac{1}{2}} B_H$ such that $g = C_\lambda \circ f - C_\lambda \circ f_{P,\lambda}$ we have*

$$\|g\|_\infty \leq 3\Big(\frac{\mathbb{E}_P g}{\lambda}\Big)^{\frac{\alpha}{2}} + \Big(\frac{a(\lambda)}{\lambda}\Big)^{\frac{\alpha}{2}} + 2 \quad \text{and}$$

$$\|f\|_H \leq \Big(\frac{a(\lambda) + \mathbb{E}_P g}{\lambda}\Big)^{1/2}.$$

*Proof.* Let us write $\varepsilon := \mathbb{E}_P g$. Then we have

$$\lambda\|f\|_H^2 \leq \lambda\|f\|_H^2 + \mathcal{R}_{L,P}(f) - \mathcal{R}_{L,P}^*$$
$$= \lambda\|f_{P,\lambda}\|^2 + \mathcal{R}_{L,P}(f_{P,\lambda}) - \mathcal{R}_{L,P}^* + \varepsilon$$
$$= a(\lambda) + \varepsilon,$$

which establishes the second assertion. Consequently, $\|L \circ f\|_\infty \leq 1 + \|f\|_\infty^\alpha$ yields

$$\|C_\lambda \circ f\|_\infty \leq \lambda\|f\|_H^2 + \|L \circ f\|_\infty \leq a(\lambda) + \varepsilon + \Big(\frac{a(\lambda)}{\lambda}\Big)^{\frac{\alpha}{2}} + \Big(\frac{\varepsilon}{\lambda}\Big)^{\frac{\alpha}{2}} + 1.$$

Analogously, we obtain $\|C_\lambda \circ f_{P,\lambda}\|_\infty \leq a(\lambda) + \big(\frac{a(\lambda)}{\lambda}\big)^{\frac{\alpha}{2}} + 1$, and therefore we find

$$\|g\|_\infty \leq \max\big(\|C_\lambda \circ f\|_\infty, \|C_\lambda \circ f_{P,\lambda}\|_\infty\big) \leq \varepsilon + \Big(\frac{a(\lambda)}{\lambda}\Big)^{\frac{\alpha}{2}} + \Big(\frac{\varepsilon}{\lambda}\Big)^{\frac{\alpha}{2}} + 2,$$

where in the last step we used $a(\lambda) \leq 1$. Now, $f \in \lambda^{-1/2} B_H$ implies that $\|f\|_\infty \leq \lambda^{-1/2}$ and an easy calculation shows that $2 + \lambda^{-\alpha/2} \leq 3\lambda^{-\frac{\alpha}{2-\alpha}}$. Therefore we obtain

$$\varepsilon \leq \mathbb{E}_P C_\lambda \circ f = \lambda\|f\|_H^2 + \mathcal{R}_{L,P}(f) \leq 2 + \|f\|_\infty^\alpha \leq 2 + \lambda^{-\frac{\alpha}{2}} \leq 3\lambda^{-\frac{\alpha}{2-\alpha}}.$$

From this we easily obtain $\varepsilon \leq 3^{1-\frac{\alpha}{2}}(\frac{\varepsilon}{\lambda})^{\frac{\alpha}{2}} \leq 2(\frac{\varepsilon}{\lambda})^{\frac{\alpha}{2}}$, which gives the assertion. $\square$

We now prove that a variance bound of the form (5) assumed in Theorem 2.1 implies a variance bound of the form (7) assumed in Theorem 3.1:

**Lemma 4.2.** *Let $P$ be a distribution on $X \times Y$ and suppose that there exist constants $v \geq 0$, $c \geq 1$, and $\vartheta \in [0,1]$ such that the variance bound assumption (5) is satisfied for some $0 < \lambda < 1$ and all $f \in \lambda^{-\frac{1}{2}} B_H$. Then for all $g \in \mathcal{G}(\lambda)$ we have*

$$\mathbb{E}_P g^2 \leq 16c\bigg(\Big(\frac{\mathbb{E}_P g}{\lambda}\Big)^{\frac{1}{2}} + \Big(\frac{a(\lambda)}{\lambda}\Big)^{1/2} + 1\bigg)^v \Big((\mathbb{E}_P g)^\vartheta + 2a^\vartheta(\lambda)\Big).$$

*Proof.* We use the shorthand notation $\mathbb{E}$ for $\mathbb{E}_P$. For $g \in \mathcal{G}(\lambda)$ pick an $f \in \lambda^{-\frac{1}{2}} B_H$ such that $g = C_\lambda \circ f - C_\lambda \circ f_{P,\lambda}$. Now observe that

$$\mathbb{E} g^2 = \mathbb{E}\big(C_\lambda \circ f - C_\lambda \circ f_{P,\lambda}\big)^2$$
$$= \mathbb{E}\big(\lambda\|f\|^2 - \lambda\|f_{P,\lambda}\|^2 + L \circ f - L \circ f_{P,\lambda}\big)^2$$
$$\leq 2\mathbb{E}\big(\lambda\|f\|^2 - \lambda\|f_{P,\lambda}\|^2\big)^2 + 2\mathbb{E}\big(L \circ f - L \circ f_{P,\lambda}\big)^2$$
$$\leq 2\lambda^2\|f\|^4 + 2\lambda^2\|f_{P,\lambda}\|^4 + 2\mathbb{E}\big(L \circ f - L \circ f_{P,\lambda}\big)^2$$
$$\leq 4\mathbb{E}\big(L \circ f - L \circ f_{L,P}^*\big)^2 + 4\mathbb{E}\big(L \circ f_{L,P}^* - L \circ f_{P,\lambda}\big)^2 + 2\lambda^2\|f\|^4 + 2\lambda^2\|f_{P,\lambda}\|^4.$$

Denote $C := \max\left(\|f\|_\infty + 1, \|f_{P,\lambda}\|_\infty + 1\right)$. Then the assumption (5) and $a^\vartheta + b^\vartheta \leq 2(a+b)^\vartheta$ for all $a, b \geq 0$, imply that

$$\mathbb{E}(L \circ f - L \circ f^*_{L,P})^2 + \mathbb{E}(L \circ f^*_{L,P} - L \circ f_{P,\lambda})^2$$
$$\leq 2cC^v\left(\mathbb{E}(L \circ f - L \circ f^*_{L,P}) + \mathbb{E}(L \circ f_{P,\lambda} - L \circ f^*_{L,P})\right)^\vartheta.$$

Since $\lambda^2 \|f\|^4 \leq 1$ and $\lambda^2 \|f_{P,\lambda}\|^4 \leq 1$ we hence obtain

$$\mathbb{E}g^2 \leq 8cC^v\left(\mathbb{E}(L \circ f - L \circ f^*_{L,P}) + \mathbb{E}(L \circ f_{P,\lambda} - L \circ f^*_{L,P})\right)^\vartheta + 2\lambda^2\|f\|^4 + 2\lambda^2\|f_{P,\lambda}\|^4$$
$$\leq 8cC^v\left(\mathbb{E}(L \circ f - L \circ f^*_{L,P}) + \mathbb{E}(L \circ f_{P,\lambda} - L \circ f^*_{L,P})\right)^\vartheta + 4\left(\lambda^2\|f\|^4 + \lambda^2\|f_{P,\lambda}\|^4\right)^\vartheta$$
$$\leq 16cC^v\left(\mathbb{E}(L \circ f - L \circ f^*_{L,P}) + \mathbb{E}(L \circ f_{P,\lambda} - L \circ f^*_{L,P}) + \lambda^2\|f\|^4 + \lambda^2\|f_{P,\lambda}\|^4\right)^\vartheta$$
$$= 16cC^v\left(\mathbb{E}g + 2\mathbb{E}(L \circ f_{P,\lambda} - L \circ f^*_{L,P}) + 2\lambda\|f_{P,\lambda}\|^2\right)^\vartheta$$
$$\leq 16cC^v\left((\mathbb{E}g)^\vartheta + 2a^\vartheta(\lambda)\right).$$

What is left is to bound $C$ in the right hand side of this inequality. To that end observe that Lemma 4.1 implies

$$\|f\|_\infty \leq \|f\|_H \leq \left(\frac{a(\lambda) + \mathbb{E}g}{\lambda}\right)^{1/2}$$

and

$$\|f_{P,\lambda}\|_\infty \leq \|f_{P,\lambda}\|_H \leq \left(\frac{a(\lambda)}{\lambda}\right)^{1/2} \leq \left(\frac{a(\lambda) + \mathbb{E}g}{\lambda}\right)^{1/2}$$

so that we can bound

$$C = \max\left(\|f\|_\infty + 1, \|f_{P,\lambda}\|_\infty + 1\right)$$
$$\leq \left(\frac{a(\lambda) + \mathbb{E}g}{\lambda}\right)^{1/2} + 1 \leq \left(\frac{\mathbb{E}g}{\lambda}\right)^{1/2} + \left(\frac{a(\lambda)}{\lambda}\right)^{1/2} + 1. \quad \square$$

The following lemma relates the covering numbers of $B_H$ with $\omega_{P,n}(\mathcal{G}(\lambda), \varepsilon)$:

**Lemma 4.3.** *Let $n \in \mathbb{N}$, and assume that there are constants $a \geq 1$ and $p \in (0,2)$ such that for all $\delta > 0$, we have*

$$\sup_{T \in Z^n} \log \mathcal{N}(B_H, \delta, L_2(T_X)) \leq a\delta^{-p}.$$

*Then there is a constant $c_{L,p} > 0$ depending only on $L$ and $p$ such that for all distributions $P$ on $X \times Y$, and all $\lambda \in (0,1]$, $\varepsilon > 0$ we have*

$$\omega_{P,n}(\mathcal{G}(\lambda), \varepsilon) \leq c_{L,p} \max\left\{\left(\frac{a(\lambda) + \varepsilon}{\lambda} + 1\right)^{\frac{\alpha p}{4}} \tau_\varepsilon^{\frac{2-p}{4}} \left(\frac{a}{n}\right)^{\frac{1}{2}}, \left(\frac{a(\lambda) + \varepsilon}{\lambda} + 1\right)^{\frac{\alpha}{2}} \left(\frac{a}{n}\right)^{\frac{2}{2+p}}\right\},$$

*where $\tau_\varepsilon \geq \sup_{g \in \mathcal{G}_\varepsilon} \mathbb{E}_P g^2$ and $\mathcal{G}_\varepsilon := \{g \in \mathcal{G}(\lambda) : \mathbb{E}_P g \leq \varepsilon\}$.*

*Proof.* Our first goal is to bound the covering numbers of $\mathcal{G}_\varepsilon$. To this end recall that for $g := C_\lambda \circ f - C_\lambda \circ f_{P,\lambda} \in \mathcal{G}_\varepsilon$, Lemma 4.1 shows that $\|f\|_H \leq \left(\frac{a(\lambda)+\varepsilon}{\lambda}\right)^{1/2} =: \Lambda$.

With the help of the auxiliary sets $\hat{\mathcal{G}}_\varepsilon := \{C_\lambda \circ f : f \in \Lambda B_H\}$ and $\mathcal{H} := \{L \circ f : f \in \Lambda B_H\}$ we thus obtain

$$\log \mathcal{N}(\mathcal{G}_\varepsilon, 2\delta, L_2(T)) \leq \log \mathcal{N}(\hat{\mathcal{G}}_\varepsilon, 2\delta, L_2(T))$$
$$\leq \log\left(\frac{1}{\delta} + 1\right) + \log \mathcal{N}(\mathcal{H}, \delta, L_2(T))$$
$$\leq \log\left(\frac{1}{\delta} + 1\right) + \log \mathcal{N}\left(\Lambda B_H, \frac{\delta}{|L_{|[-\Lambda,\Lambda]}|_1}, L_2(T_X)\right).$$

Furthermore, the Lipschitz assumption (1) implies the right hand side is bounded by

$$\log\left(\frac{1}{\delta} + 1\right) + \log \mathcal{N}\left(B_H, \frac{\delta}{c_L \Lambda^\alpha}, L_2(T_X)\right).$$

Consequently, there is a constant $\tilde{c}_{L,p} > 0$ depending only on $L$ and $p$ such that for all $\delta > 0$ we have

$$\sup_{T \in Z^n} \log \mathcal{N}(\mathcal{G}_\varepsilon, \delta, L_2(T)) \leq a\,\tilde{c}_{L,p}\left(\frac{a(\lambda) + \varepsilon}{\lambda}\right)^{\frac{\alpha p}{2}} \delta^{-p} \leq a\,\tilde{c}_{L,p}\left(\frac{a(\lambda) + \varepsilon}{\lambda} + 1\right)^{\frac{\alpha p}{2}} \delta^{-p}.$$

By symmetrization, and the proofs of [9, Lem. 2.5] and [12, Prop. 5.7] we thus find

$$\omega_{P,n}(\mathcal{G}(\lambda), \varepsilon) \leq c_{L,p} \max\left\{\left(\frac{a(\lambda)+\varepsilon}{\lambda}+1\right)^{\frac{\alpha p}{4}} \tau_\varepsilon^{\frac{2-p}{4}} \left(\frac{a}{n}\right)^{\frac{1}{2}}, \left(\frac{a(\lambda)+\varepsilon}{\lambda}+1\right)^{\frac{\alpha}{2}} \left(\frac{a}{n}\right)^{\frac{2}{2+p}}\right\}.$$

□

*Proof of Theorem 2.1.* Let $g := C_\lambda \circ f - C_\lambda \circ f_{P,\lambda}$ for some $f \in \lambda^{-1/2} B_H$. Lemma 4.1 implies that we have a supremum bound

$$\|g\|_\infty \leq 3\left(\frac{\mathbb{E}_P g}{\lambda}\right)^{\frac{\alpha}{2}} + \left(\frac{a(\lambda)}{\lambda}\right)^{\frac{\alpha}{2}} + 2.$$

Because of the variance bound assumption (5), Lemma 4.2 implies we have a variance bound of the form

$$\mathbb{E}_P g^2 \leq 16c\left(\left(\frac{\mathbb{E}_P g}{\lambda}\right)^{\frac{1}{2}} + \left(\frac{a(\lambda)}{\lambda}\right)^{1/2} + 1\right)^v \left((\mathbb{E}_P g)^\vartheta + 2a^\vartheta(\lambda)\right)$$
$$\leq 48c\left(\left(\frac{\mathbb{E}_P g}{\lambda}\right)^{\frac{\alpha}{2}} + \left(\frac{a(\lambda)}{\lambda}\right)^{\frac{\alpha}{2}} + 1\right)^{\frac{v}{\alpha}} \left((\mathbb{E}_P g)^\vartheta + 2a^\vartheta(\lambda)\right)$$
$$\leq \left(3\left(\frac{\mathbb{E}_P g}{\lambda}\right)^{\frac{\alpha}{2}} + \left(\frac{a(\lambda)}{\lambda}\right)^{\frac{\alpha}{2}} + 2\right)^{\frac{v}{\alpha}} \left(48c(\mathbb{E}_P g)^\vartheta + 96ca^\vartheta(\lambda)\right).$$

Therefore we have variance and supremum bounds of the form (7) and (6) with the values $b = 3\lambda^{-\frac{\alpha}{2}}$, $\beta = \frac{\alpha}{2}$, $B = \left(\frac{a(\lambda)}{\lambda}\right)^{\frac{\alpha}{2}} + 2$, $w = 48c$, $\nu = \frac{v}{\alpha}$, and $W = 96ca^\vartheta(\lambda)$. Denote $\tau_\varepsilon := 3^4 2^6 c\lambda^\vartheta \left(\frac{a(\lambda)+\epsilon}{\lambda}+1\right)^{\vartheta+\frac{v}{2}}$. Then for $g \in \mathcal{G}(\lambda)$ with $\mathbb{E}_P g \leq \varepsilon$ we obtain

$$\mathbb{E}_P g^2 \leq (b\varepsilon^\beta + B)^\nu (w\varepsilon^\vartheta + W)$$
$$= 48c\left(3\left(\frac{\varepsilon}{\lambda}\right)^{\frac{\alpha}{2}} + \left(\frac{a(\lambda)}{\lambda}\right)^{\frac{\alpha}{2}} + 2\right)^{\frac{v}{\alpha}} \left(\varepsilon^\vartheta + 2a^\vartheta(\lambda)\right)$$
$$\leq 96 \cdot 9c\left(\left(\frac{\varepsilon}{\lambda}\right)^{\frac{\alpha}{2}} + \left(\frac{a(\lambda)}{\lambda}\right)^{\frac{\alpha}{2}} + 1\right)^{\frac{v}{\alpha}} \left(\varepsilon^\vartheta + a^\vartheta(\lambda)\right)$$
$$\leq 96 \cdot 9 \cdot 3 \cdot 2c\left(\frac{a(\lambda)+\varepsilon}{\lambda}+1\right)^{\frac{v}{2}} \left(a(\lambda)+\varepsilon\right)^\vartheta$$
$$\leq \tau_\varepsilon.$$

Consequently we can apply Lemma 4.3 to obtain that

$$\omega_{P,n}(\mathcal{G}(\lambda),\varepsilon)$$
$$\leq c_{L,p} \max\left\{\left(\frac{a(\lambda)+\varepsilon}{\lambda}+1\right)^{\frac{p}{4}} \tau_\varepsilon^{\frac{2-p}{4}} \left(\frac{a}{n}\right)^{\frac{1}{2}}, \left(\frac{a(\lambda)+\varepsilon}{\lambda}+1\right)^{\frac{1}{2}} \left(\frac{a}{n}\right)^{\frac{2}{2+p}}\right\}$$
$$\leq c_{L,p}\sqrt{c} \max\left\{\lambda^{\frac{\vartheta(2-p)}{4}} \left(\frac{a(\lambda)+\varepsilon}{\lambda}+1\right)^{\frac{2\alpha p+(v+2\vartheta)(2-p)}{8}} \left(\frac{a}{n}\right)^{\frac{1}{2}}, \left(\frac{a(\lambda)+\varepsilon}{\lambda}+1\right)^{\frac{\alpha}{2}} \left(\frac{a}{n}\right)^{\frac{2}{2+p}}\right\}.$$

We also bound the terms

$$\sqrt{\frac{2x(b\varepsilon^\beta + B)^\nu(w\varepsilon^\vartheta + W)}{n}} \leq \sqrt{\frac{2x\tau_\varepsilon}{n}} = 384\sqrt{2c}\lambda^{\frac{\vartheta}{2}} \left(\frac{a(\lambda)+\varepsilon}{\lambda}+1\right)^{\frac{\vartheta}{2}+\frac{v}{4}} \sqrt{\frac{x}{n}}$$

and

$$\frac{2x(b\varepsilon^\beta + B)}{n} = \frac{2x}{n}\left(3\left(\frac{\epsilon}{\lambda}\right)^{\frac{\alpha}{2}} + \left(\frac{a(\lambda)}{\lambda}\right)^{\frac{\alpha}{2}} + 2\right) \leq \frac{24x}{n}\left(\frac{a(\lambda)+\epsilon}{\lambda}+1\right)^{\frac{\alpha}{2}}$$

and then observe that Theorem 3.1 implies that there is a constant $K \geq 1$ such that

$$\Pr{}^*\Big(T \in (X \times Y)^n : \mathcal{R}_{C_\lambda,P}(f_{T,\lambda}) < \mathcal{R}_{C_\lambda,P}(f_{P,\lambda}) + \tilde{\varepsilon}\Big) \geq 1 - e^{-x},$$

whenever

$$\tilde{\varepsilon} \geq K\max\left\{\lambda^{\frac{\vartheta(2-p)}{4}} \left(\frac{a(\lambda)+\tilde{\varepsilon}}{\lambda}+1\right)^{\frac{2\alpha p+(v+2\vartheta)(2-p)}{8}} \left(\frac{a}{n}\right)^{\frac{1}{2}}, \left(\frac{a(\lambda)+\tilde{\varepsilon}}{\lambda}+1\right)^{\frac{\alpha}{2}} \left(\frac{a}{n}\right)^{\frac{2}{2+p}},\right.$$
$$\left. \lambda^{\frac{\vartheta}{2}} \left(\frac{a(\lambda)+\tilde{\varepsilon}}{\lambda}+1\right)^{\frac{\vartheta}{2}+\frac{v}{4}} \left(\frac{x}{n}\right)^{\frac{1}{2}}, \left(\frac{a(\lambda)+\tilde{\varepsilon}}{\lambda}+1\right)^{\frac{\alpha}{2}} \frac{x}{n}\right\}.$$

If we further constrain by $\tilde{\varepsilon} \geq a(\lambda) + \lambda$ we find that it is sufficient to satisfy

$$\tilde{\varepsilon} \geq \max\left\{a(\lambda)+\lambda, K\lambda^{\frac{\vartheta(2-p)}{4}} \left(\frac{\tilde{\varepsilon}}{\lambda}\right)^{\frac{2\alpha p+(v+2\vartheta)(2-p)}{8}} \left(\frac{a}{n}\right)^{\frac{1}{2}}, K\left(\frac{\tilde{\varepsilon}}{\lambda}\right)^{\frac{\alpha}{2}} \left(\frac{a}{n}\right)^{\frac{2}{2+p}},\right.$$
$$\left. K\lambda^{\frac{\vartheta}{2}} \left(\frac{\tilde{\varepsilon}}{\lambda}\right)^{\frac{\vartheta}{2}+\frac{v}{4}} \left(\frac{x}{n}\right)^{\frac{1}{2}}, K\left(\frac{\tilde{\varepsilon}}{\lambda}\right)^{\frac{\alpha}{2}} \frac{x}{n}\right\}.$$

Since $\vartheta \in (0,1]$ and $v \in [0,2]$ it follows that $0 < v + 2\vartheta \leq 4$ which implies that $\frac{2\alpha p+(v+2\vartheta)(2-p)}{8} \leq 1$ and $\frac{\vartheta}{2} + \frac{v}{4} \leq 1$. Consequently we find that it is sufficient to satisfy

$$\tilde{\varepsilon} \geq \max\left\{a(\lambda)+\lambda, \left(\frac{K^2 a}{\lambda^{\frac{2\alpha p+v(2-p)}{4}}n}\right)^{\frac{4}{8-2\alpha p-(v+2\vartheta)(2-p)}}, \left(\frac{K^{\frac{2+p}{2}}a}{\lambda^{\frac{\alpha(2+p)}{4}}n}\right)^{\frac{4}{(2+p)(2-\alpha)}},\right.$$
$$\left.\left(\frac{K^2 x}{\lambda^{\frac{v}{2}}n}\right)^{\frac{2}{4-(v+2\vartheta)}}, \left(\frac{Kx}{\lambda^{\frac{\alpha}{2}}n}\right)^{\frac{2}{2-\alpha}}\right\}.$$

Therefore we find that (with a change in the value of the constant $K$) if

$$\varepsilon \geq \max\left\{a(\lambda)+\lambda, \left(\frac{Ka}{\lambda^{\frac{2\alpha p+v(2-p)}{4}}n}\right)^{\frac{4}{8-2\alpha p-(v+2\vartheta)(2-p)}}, \left(\frac{Ka}{\lambda^{\frac{\alpha(2+p)}{4}}n}\right)^{\frac{4}{(2+p)(2-\alpha)}},\right.$$
$$\left.\left(\frac{Kx}{\lambda^{\frac{v}{2}}n}\right)^{\frac{2}{4-(v+2\vartheta)}}, \left(\frac{Kx}{\lambda^{\frac{\alpha}{2}}n}\right)^{\frac{2}{2-\alpha}}\right\}$$

then

$$\mathcal{R}_{L,P}(f_{T,\lambda}) \leq \lambda \|f_{T,\lambda}\|_H^2 + \mathcal{R}_{L,P}(f_{T,\lambda}) = \mathcal{R}_{C_\lambda,P}(f_{T,\lambda}) < \mathcal{R}_{C_\lambda,P}(f_{P,\lambda}) + \varepsilon$$
$$= a(\lambda) + \mathcal{R}_{L,P}^* + \varepsilon$$

holds with probability not less than $1 - e^{-x}$. □

## 5. Examples

Here we perform the analysis mentioned in Examples 2.2 and 2.3. Let us first apply the oracle inequality to bound $\mathcal{R}_{L_\alpha,P}(f_{T,\lambda}) - \mathcal{R}_{L_\alpha,P}^*$ with high probability. To that end we now derive some variance bounds. First observe that [11, Table 3] shows that the modulus of convexity $\delta_{\psi_\alpha|[-B,B]}(\varepsilon)$ of the function $\psi_\alpha : t \mapsto |t|^\alpha$ restricted to the interval $[-B, B]$ satisfies

$$(8) \qquad \delta_{\psi_\alpha|[-B,B]}(\varepsilon) \geq \frac{\alpha(\alpha-1)}{8} B^{\alpha-2} \varepsilon^2$$

Consequently [2, Lemma 15] implies that modulus of convexity of $\mathcal{R}_{L_\alpha,P}$ for functions satisfying $\|f\|_\infty \leq B$ is bounded below by $\frac{\alpha(\alpha-1)}{8} 2^{\alpha-2} B^{\alpha-2} \varepsilon^2 \geq \frac{\alpha(\alpha-1)}{16} \times B^{\alpha-2} \varepsilon^2$. Moreover, the mean value theorem implies that

$$\left||t_1 - y|^\alpha - |t_2 - y|^\alpha\right| \leq \alpha \Big(\max(t_1 + 1, t_2 + 1)\Big)^{\alpha-1} |t_1 - t_2|$$

so that the loss function $f \mapsto L_\alpha(y, f(x))$ has a Lipschitz constant less than $\alpha \big(\max\{\|f_1\|_\infty, \|f_2\|_\infty\} + 1\big)^{\alpha-1}$. Now let

$$f_{L_\alpha,P}^* \in \arg\min\{\mathcal{R}_{L_\alpha,P}(f) | f : X \to \mathbb{R} \text{ measurable}\}$$

and define $g_f(x, y) := |f(x) - y|^\alpha - |f_{L_\alpha,P}^*(x) - y|^\alpha$. Then the extension mentioned after the statement of [2, Lemma 14] to non-margin loss functions implies that we have the variance bound

$$\mathbb{E}g_f^2 \leq \frac{8\alpha}{(\alpha-1)} \cdot \frac{\big(\max\{\|f\|_\infty, \|f_{L_\alpha,P}^*\|_\infty\} + 1\big)^{2\alpha-2}}{\big(\max\{\|f\|_\infty, \|f_{L_\alpha,P}^*\|_\infty\}\big)^{\alpha-2}} \mathbb{E}g_f$$
$$\leq \frac{8\alpha}{(\alpha-1)} \big(\max\{\|f\|_\infty, \|f_{L_\alpha,P}^*\|_\infty\} + 1\big)^\alpha \mathbb{E}g_f..$$

Observe that the right hand side of these bounds goes to $\infty$ as $\alpha \to 1$ since $\psi_1$ is not strictly convex. Also note that such a bound, but with different constants, follows directly from [11, Equation 28]. Since $\|f_{L_\alpha,P}^*\|_\infty \leq 1$ we then obtain

$$\mathbb{E}g_f^2 \leq \frac{8\alpha}{(\alpha-1)} \big(\|f\|_\infty + 2\big)^\alpha \mathbb{E}g_f..$$

Therefore we can apply Theorem 2.1 with $v = \alpha$ and $\vartheta = 1$ to obtain that there exists a constant $K_\alpha \geq 1$ such that for all $0 < \lambda \leq 1$, $\varepsilon > 0$, $x \geq 1$ satisfying

$$(9) \qquad \varepsilon \geq \max\left\{a_\alpha(\lambda) + \lambda, \Big(\frac{K_\alpha a}{\lambda^{\frac{\alpha(2+p)}{4}} n}\Big)^{\frac{4}{(2+p)(2-\alpha)}}, \Big(\frac{K_\alpha x}{\lambda^{\frac{\alpha}{2}} n}\Big)^{\frac{2}{2-\alpha}}\right\},$$

we have

$$\text{(10)} \qquad \Pr{}^*\left(T \in Z^n : \mathcal{R}_{L_\alpha,P}(f_{T,\lambda}) - \mathcal{R}^*_{L_\alpha,P} < a_\alpha(\lambda) + \varepsilon\right) \geq 1 - e^{-x}$$

where $a_\alpha(\cdot)$ is the approximation error function defined with respect to the risk $\mathcal{R}_{L_\alpha,P}$.

In [13] it was shown that the assumption $f^*_{L_\alpha,P} \in H$ implies that $a_\alpha(\lambda) \leq \lambda \|f^*_{L_\alpha,P}\|^2_H$ for all $\lambda > 0$. We assume without loss of generality that $a_\alpha(\lambda) \leq \lambda$. Let us first consider when $\alpha = 2$. If we now assume $(\lambda_n)$ is a strictly positive null-sequence with $\lambda_n^{1+p/2} n \to \infty$ then it is easy from the convention (3) applied to the inequality (9) that our learning rate is of the form $\lambda_n$ thus finishing the proof for Example 2.2. Now consider the case $1 < \alpha < 2$. Then (9) becomes

$$\text{(11)} \qquad \varepsilon \geq \max\left\{a_\alpha(\lambda) + \lambda, \lambda^{-\frac{\alpha}{2-\alpha}}\left(\frac{K_\alpha a}{n}\right)^{\frac{4}{(2+p)(2-\alpha)}}, \lambda^{-\frac{\alpha}{2-\alpha}}\left(\frac{K_\alpha x}{n}\right)^{\frac{2}{2-\alpha}}\right\}.$$

Moreover when $n \geq K_\alpha a$ elementary calculations show that it is sufficient to satisfy

$$\text{(12)} \qquad \varepsilon \geq a_\alpha(\lambda) + \lambda + \lambda^{-\frac{\alpha}{2-\alpha}} x^{\frac{2}{2-\alpha}} \left(\frac{K_\alpha a}{n}\right)^{\frac{4}{(2+p)(2-\alpha)}}..$$

If we now assume $\lambda = n^{-\kappa}$. Then elementary calculations show that we obtain the rate $n^{-\kappa}$ independently of the value $\alpha$ when $\kappa \leq \frac{2}{2+p}$ and when $\kappa > \frac{2}{2+p}$ we obtain the rate $n^{-\frac{2}{2+p} + \frac{\alpha}{2-\alpha}(\kappa - \frac{2}{2+p})}$.

Let us now assume that the conditional distributions $P(y|x)$ are symmetric. We now proceed to derive a calibration inequality

$$\mathcal{R}_{L_2,P}(f_{T,\lambda}) - \mathcal{R}^*_{L_2,P} \leq \Psi(\mathcal{R}_{L_\alpha,P}(f_{T,\lambda}) - \mathcal{R}^*_{L_\alpha,P})$$

so that we can apply the bounds on $\mathcal{R}_{L_\alpha,P}(f_{T,\lambda}) - \mathcal{R}^*_{L_\alpha,P})$ defined by (10) and (12) to obtain bounds on $\mathcal{R}_{L_2,P}(f_{T,\lambda}) - \mathcal{R}^*_{L_2,P}$ in terms of $\alpha$. Since we will need results and notations from [11] we first give a brief outline of its content. Consider a loss function $L$ and a measure $Q$ on $Y$. Then the associated inner risk is defined as

$$\mathcal{C}_{L,Q}(t) = \int_Y L(y,t) dQ(y), \quad t \in \mathbb{R},$$

and can be used to compute the risk

$$\mathcal{R}_{L,P}(f) = \int_X \mathcal{C}_{L,P(\cdot|x)}(f(x)) dP_X(x).$$

The minimal inner risk is defined as $\mathcal{C}^*_{L,Q} := \inf_{t \in \mathbb{R}} \mathcal{C}_{L,Q}(t)$. Consider now another loss function $\acute{L}$. Then the calibration function $\delta_{\max,L,\acute{L}}(\varepsilon,Q)$ is defined as the largest function comparing the excess inner risks, i.e.

$$\delta_{\max,L,\acute{L}}(\mathcal{C}_{L,Q}(t) - \mathcal{C}^*_{L,Q}, Q) \leq \mathcal{C}_{\acute{L},Q}(t) - \mathcal{C}^*_{\acute{L},Q}.$$

We shall also find it convenient to consider the template loss $L_{\text{mean}}$ introduced in [11] and defined by

$$L_{\text{mean}}(Q,t) := |\mathbb{E}Q - t|, \quad t \in \mathbb{R}$$

and its inner risk

$$\mathcal{C}_{L_{\text{mean}},Q}(t) = \int_Y |\mathbb{E}Q - t| dQ(y), \quad t \in \mathbb{R}.$$

We can now proceed to derive the appropriate calibration inequality function $\Psi$ for comparing $L_2$ and $L_\alpha$. Since $P(y|x)$ is symmetric for all $x$, [11, Theorem 3.23] implies that we have mean calibration with calibration function bounded below by

$$\delta_{\max,L_{\mathrm{mean}},L_\alpha}(\varepsilon,Q) \geq \delta_{\psi_\alpha|[-(2+\varepsilon),2+\varepsilon]}(2\varepsilon)$$

where $\delta_{\psi_\alpha|[-(2+\varepsilon),2+\varepsilon]}$ is the modulus of convexity of the function $\psi_\alpha$ restricted to the interval $[-(2+\varepsilon), 2+\varepsilon]$. By (8) we then obtain

$$\delta_{\max,L_{\mathrm{mean}},\psi_\alpha}(\varepsilon,Q) \geq \frac{\alpha(\alpha-1)}{2}(2+\varepsilon)^{\alpha-2}\varepsilon^2.$$

Since [11, Equation (38)] states $\delta_{\max,L_2,L_\alpha}(\varepsilon,Q) = \delta_{\max,L_{\mathrm{mean}},L_\alpha}(\sqrt{\varepsilon},Q)$ we find

$$\delta_{\max,L_2,L_\alpha}(\varepsilon,Q) \geq \frac{\alpha(\alpha-1)}{2}(2+\sqrt{\varepsilon})^{\alpha-2}\varepsilon.$$

We now seek to apply [11, Theorem 2.13]. In that notation we bound

$$B_f = \sup_x |f(x) - \mathbb{E}(y|x)|^2 \leq |\|f\|_\infty + 1|^2.$$

Denote $\phi(\varepsilon) := \frac{\alpha(\alpha-1)}{2}(2+\sqrt{\varepsilon})^{\alpha-2}\varepsilon$. Then since

$$\frac{d}{d\varepsilon}\left((2+\sqrt{\varepsilon})^{\alpha-2}\varepsilon\right) = (2+\sqrt{\varepsilon})^{\alpha-3}(2+\frac{\alpha}{2}\sqrt{\varepsilon}) > 0$$

and

$$\frac{d^2}{d\varepsilon^2}\left((2+\sqrt{\varepsilon})^{\alpha-2}\varepsilon\right) = (\alpha-2)\varepsilon^{-\frac{1}{2}}\left(\frac{3}{2}+\frac{\alpha}{4}\sqrt{\varepsilon}\right)(2+\sqrt{\varepsilon})^{\alpha-4} \leq 0$$

we conclude that $\phi$ is strictly monotonically increasing and concave. It follows that

$$\phi^{**}_{B_f}(\varepsilon) \geq \phi^{**}_{|\|f\|_\infty+1|^2}(\varepsilon) = \frac{\phi(|\|f\|_\infty+1|^2)}{|\|f\|_\infty+1|^2}\varepsilon = \frac{\alpha(\alpha-1)}{2}(3+\|f\|_\infty)^{\alpha-2}\varepsilon$$

where ** denotes the Fenchel-Legendre bi-conjugate operation (see e.g. [10]). It then follows from [11, Theorem 2.13] that

$$(13) \quad \mathcal{R}_{L_2,P}(f) - \mathcal{R}^*_{L_2,P} \leq \frac{2}{\alpha(\alpha-1)}(3+\|f\|_\infty)^{2-\alpha}(\mathcal{R}_{L_\alpha,P}(f) - \mathcal{R}^*_{L_\alpha,P})$$

for all bounded measurable functions $f$. Note that the constant in this inequality goes to $\infty$ as $\alpha$ goes to 1. The deeper reason for this behaviour is that $\psi_\alpha$ is strictly convex when $\alpha > 1$ but not strictly convex when $\alpha = 1$ as discussed in [11].

We conclude from inequalities (13) and (10) that whenever (12) is satisfied that with probability greater than $1 - e^{-x}$ we have

$$\mathcal{R}_{L_2,P}(f_{T,\lambda}) - \mathcal{R}^*_{L_2,P} \leq \frac{4}{\alpha(\alpha-1)}(3+\|f_{T,\lambda}\|_\infty)^{2-\alpha}\varepsilon.$$

However we also know from the last line of the proof of Theorem 2.1 that whenever (12) is satisfied that with probability greater than $1 - e^{-x}$

$$\|f_{T,\lambda}\|_\infty \leq \|f_{T,\lambda}\|_H \leq \sqrt{\frac{a_\alpha(\lambda)+\varepsilon}{\lambda}} \leq \sqrt{2}\sqrt{\frac{\varepsilon}{\lambda}}.$$

Now since $1 \leq \frac{\varepsilon}{\lambda}$ when (12) is satisfied it follows that

$$(3 + \|f_{T,\lambda}\|_\infty)^{2-\alpha} \leq \left(3 + \sqrt{2}\sqrt{\frac{\varepsilon}{\lambda}}\right)^{2-\alpha} \leq \left(3 + \sqrt{2}\right)\left(\frac{\varepsilon}{\lambda}\right)^{1-\frac{\alpha}{2}}$$

so that with probability greater than $1 - 2e^{-x}$ we have

$$\mathcal{R}_{L_2,P}(f_{T,\lambda}) - \mathcal{R}^*_{L_2,P} \leq \frac{4}{\alpha(\alpha-1)}(3+\sqrt{2})\lambda^{\frac{\alpha}{2}-1}\varepsilon^{2-\frac{\alpha}{2}}.$$

If we now apply the inequality $a_\alpha(\lambda) \leq \lambda$ and let

$$\varepsilon := 2\lambda + \lambda^{-\frac{\alpha}{2-\alpha}}x^{\frac{2}{2-\alpha}}\left(\frac{K_\alpha a}{n}\right)^{\frac{4}{(2+p)(2-\alpha)}},$$

then we see that with probability greater than $1 - 2e^{-x}$ we have

$$\mathcal{R}_{L_2,P}(f_{T,\lambda}) - \mathcal{R}^*_{L_2,P}$$
$$\leq \frac{4}{\alpha(\alpha-1)}(3+\sqrt{2})\lambda^{\frac{\alpha}{2}-1}\left(2\lambda + \lambda^{-\frac{\alpha}{2-\alpha}}x^{\frac{2}{2-\alpha}}\left(\frac{K_\alpha a}{n}\right)^{\frac{4}{(2+p)(2-\alpha)}}\right)^{2-\frac{\alpha}{2}}$$
$$\leq c_\alpha\left(\lambda + \lambda^{-\frac{2}{2-\alpha}}x^{\frac{4-\alpha}{2-\alpha}}\left(\frac{K_\alpha a}{n}\right)^{\frac{2}{2+p}\frac{4-\alpha}{2-\alpha}}\right)$$

for some constant $c_\alpha$ which depends only on $\alpha$.

Now let us consider the case when $\lambda = n^{-\kappa}$. Then disregarding the constants the righthand side becomes

$$n^{-\kappa} + n^{(\kappa - \frac{2}{2+p})\frac{2}{2-\alpha} - \frac{2}{2+p}}$$

so that we obtain performance bounds of the form $n^{-\rho}$ with

$$\rho = \min\left(\kappa, \frac{2}{2+p} + \left(\frac{2}{2+p} - \kappa\right)\frac{2}{2-\alpha}\right).$$

Simple calculations show that when $\kappa \leq \frac{2}{2+p}$ then $\rho = \kappa$ independently of the value of $\alpha$ and when $\kappa > \frac{2}{2+p}$ then $\rho = \frac{2}{2+p} + (\frac{2}{2+p} - \kappa)\frac{2}{2-\alpha}$. In the latter case it is important to observe that the rates get worse as $\alpha$ increases towards 2. Indeed one can show that $\rho \leq 0$ in the interval

$$2 - (\kappa - \frac{2}{2+p})(2+p) \leq \alpha \leq 2.$$

Moreover one can see that smaller $\alpha$ minimizes the sensitivity to the degree to which $\kappa$ is greater than $\frac{2}{2+p}$.

### Acknowledgment

Ingo Steinwart thanks T. Zhang and P. Bartlett for the fruitful discussions in Eindhoven and Albuquerque that helped to develop the double localization argument. He further thanks A. Caponnetto and V. Temlyakov for the discussions on optimal rates for Sobolev spaces at the Toyota Technological Institute at Chicago.